\documentclass[12pt]{article}
\usepackage{epsfig,amsmath,amssymb,multirow}
\usepackage{authblk}
\title{A fast ADI algorithm for nonlinear Poisson equation in heterogeneous dielectric media}
\author{ Wufeng Tian\thanks{wufeng.tian@uwc.edu}}
\affil{Department of Mathematics, University of Wisconsin Colleges, WI}

\begin{document}

\maketitle

\begin{abstract}

Recently, a nonlinear Poisson equation has been introduced to model nonlinear and nonlocal hyperpolarization effects in electrostatic solute-solvent interaction for biomolecular solvation analysis. Due to a strong nonlinearity associated with the heterogeneous dielectric media, this Poisson model is difficult to solve numerically, particularly for large protein systems. A new pseudo-transient continuation approach is proposed in this paper to efficiently and stably solve the nonlinear Poisson equation. A Douglas type alternating direction implicit (ADI)  method is developed for solving the pseudo-time dependent Poisson equation. Different approximations to the dielectric profile in heterogeneous media are considered in the standard finite difference discretization. The proposed ADI scheme is validated by considering benchmark examples with exact solutions and by solvation analysis of real biomolecules with various sizes. Numerical results are in good agreement with the theoretical prediction, experimental measurements, and those obtained from the boundary value problem approach. Since the time stability of the proposed ADI scheme can be maintained even using very large time increments, it is efficient for electrostatic analysis involving hyperpolarization effects.



\end{abstract}

\section{Introduction}

With the development of theoretical methods and computational techniques over the past few 
decades, molecular modeling and simulation have become an effective and practical approach to 
mimic the behavior of molecules with biological significance. The modeling and computation of molecular surface and solute-solvent interactions, electrostatic 
interactions are of great importance in quantitative
studies of macromolecules, including proteins, DNAs, molecular motors, and viruses.
The structure, function, dynamics, and transport of macromolecules
depend on the features of molecular surfaces. Under physiological conditions,
solute-solvent interactions naturally occur in an aqueous environment. Electrostatic 
interactions are ubiquitous in nature and fundamental for chemical, biological and material sciences.

In order to practically implement a quantitative description and analysis of various important 
biological processes at the atomic level, such as DNA recognition, transcription,
translation, protein folding and protein ligand binding, the analysis of the biomolecule solvation is essential
 \cite{Baker05,Dong08,Feig04}. Biologically, solvation analysis concerns the the interaction between solute
 macromolecules and surrounding solvent molecules or ions. Mathematically, solute-solvent interactions can be represented
 via solvation energies, which contain polar and nonpolar components. The polar component, which we consider as polar
 solvation energy, is due to electrostatic interaction. Electrostatic 
interactions are universal for any system of charged or polar molecules.
The nonpolar portion consists of  surface tension effects, mechanical work and 
attractive solvent-solute dispersion interactions. The implicit solvent theory that treats solvent as a continuous dielectric medium, is one of the most commonly used approaches for studying biomolecular solvation\cite{Baker05, Roux, Sharp}.

Differential geometry based multi-scale implicit solvent models have been constructed to determine
 equilibrium properties of solvation \cite{ChenZhan1, ChenZhan2, Wei10}. These solvation models have
 implemented both Eulerian
 and Lagrangian formulations\cite{ChenZhan1, ChenZhan2}. A free energy minimization or optimization process has been
 conducted in these models based on the fundamental laws of physics. The total free energy functional for the solvation
 analysis consists of  polar component-electrostatic potential and nonpolar component-geometric effect of solvent-solute interface,
 the mechanical work of the system and the dispersive solvent-solute interaction \cite{ChenZhan1, ChenZhan2, Wei10}. 
By applying Euler-Lagrange variational approach, two coupled nonlinear partial differential equations are 
derived as governing equations\cite{ChenZhan1, ChenZhan2}. One is generalized nonlinear Poisson Boltzmann equation(NPB) for electrostatic potential 
and the other is a generalized Laplace-Beltrami equation defining the solvent-solute interface. These differential geometry based multi-scale models have been widely applied into more complex chemical and biological systems, such as molecular motors, ion channels, DNA packing and virus evolution\cite{ChenZhan1, ChenZhan2, Wei10, Crowley}.  In Poisson equation, derived by using Gauss's
Law and linear polarization, which provides a relatively simple and accurate but very much less expensive 
description of electrostatic interactions for a given charge source density. However, the Poisson equation is strictly 
valid for linear polarization in an isotropic and homogeneous dielectric medium. Unfortunately, most work based on classic Poisson theory neglects hyperpolarization, 
which is important under a strong electric field, or in highly charged nonlinear inhomogeneous media. 

A pseudo transient continuation model for the theoretical modeling of the
 biomolecular surface and solvation process has been proposed by Zhao \cite{Zhao14, Geng13}, which  completely suppressed the nonlinear instablility
through the introduction of a pseudo-time iteration process, giving rise to a time dependent NPB equation and a time dependent
Laplace-Beltrami or geometric flow equation. Unconditionally stable locally one-dimensional (LOD) methods have been developed for solving the PB equation in a one-dimensional manner\cite{Wilson16}.  By using very large time step for steady state computations, solving 1D matrices  in the implicit time stepping, these LOD pseudo-time methods are very efficient for large protein systems.
 By treating the solution of a nonlinear boundary value system as the steady state solution of a time dependent process, the over all model coupling is accomplished by the explicit Euler time integration and controlled by time increments. Compared with a relaxation
 based iterative procedure used in\cite{ChenZhan1, ChenZhan2}, this coupling is simpler and has fewer controlling
 parameters. The NPB equation can be treated in the same manner as the linearized PB equation, which overcomes the difficulties of
 conventional coupling. Computationally, the smooth solvent-solute interface possesses good differentiability and generates a
 smooth dielectric profile for the electrostatic potential so the central finite difference scheme can be employed for 
the NPB generalized equation for the discretization in space.  A fast and accurate algorithm has been proposed to solve
the time dependent Laplace-Beltrami or potential driven mean curvature flow equation in molecular surface generation\cite{Tian14}. The time splitting alternative direction implicit(ADI) method has been proposed for solving the generalized NPB equation in the pseudo-transient solvation model\cite{Zhao11, Zhao14, Geng13}.

Recently, A nonlinear Poisson equation with a nonlinear functional for a dielectric has been introduced due to 
nonlinear effects in the vicinity of the solute boundary and their impact on heterogeneous media. A new electrostatic solvation
free-energy functional partially accounting for the effect of hyperpolarizations has been proposed in electrostatic solvation analysis and treated and solved as a BVP problem by algebraic iterative process.
\cite{Hu12}.  Due to a strong nonlinearity associated with the heterogeneous dielectric media, this Poisson model is difficult to solve numerically, particularly for large protein systems. In the BVP approach \cite{Hu12}, one calculates the electrostatic potential $\phi$ and dielectric function $\epsilon$ alternatively. The computation will be regarded to be convergent, when the relative difference in the electrostatic solvation free energy of two successive steps is less than a tolerance. However, for extra large protein systems, the algebraic iterative solution process could be time consuming. Given a similar form as the Laplace-Beltrami operator in a nonlinear Poisson equation, and the efficient properties of ADI schemes, our goal is to apply the ADI method to effciently solve the nonlinear Poisson (NP) equation, which is a widely accepted model for electrostatic analysis.

The rest of paper is organized as follows: Section 2 introduces the theory and formulation of the
 the governing equation of nonlinear Poisson equation and the previous iterative numberical approach. Section 3 is devoted to the theory and formulation of the
ADI time integration.  In section 4, numerical validations based
on various benchmark biological systems are considered. Both stability and accuracy
will be examined for the proposed ADI algorithm  in solving the
Nonlinear possion equation. In section 5, we first validate our ADI scheme through the physical quantities, such as surface areas, surface enclosed volumes, electrostatic solvation energy of one atom of unit van der Waals radius and unit charge in the solvent.  Then, we apply the proposed ADI algorithm to various
atomic systems, amino acids, and protein molecules. Finally, this paper ends with
concluding remarks.

\section{Mathematical models and existing algorithms}

\subsection{Governing equation}

Consider a solute macromolecule in space surrounded by a solvent aqueous solution. To account for the nonlocal effect of hyperpolarizations in solvation analysis, a nonlinear Poisson equation has been introduced in \cite{Hu12} as a correction to the classical Poisson model
\begin{equation}\label{NPE}
-\nabla \cdot (\epsilon (|\nabla \phi({\bf r})|) \nabla \phi({\bf r}))
 = \rho_m({\bf r}),
\end{equation}
where $\phi$ is the electrostatic potential and the source term $\rho_m$ is defined as
\begin{equation}\label{RhoEq}
\rho_m({\bf r}) = \frac{4\pi e^2_c}{k_B T} \displaystyle\sum\limits_{j=1}^{N_m} q_j \delta ({\bf r} - {\bf r}_j).
\end{equation} 
Here $T$ is the temperatur, $k_B$ is the Boltzmann constant, $e_c$ is the fundamental charge, and $q_j$, in the same units as $e_c$, is the partial charge on the $j$th atom of the solute macromolecule locate at position ${\bf r}_j$.

Two general forms have been suggested in \cite{Hu12} for modeling the nonlinear dielectric function $\epsilon (|\nabla \phi|)$. One form is given as 
\begin{equation}\label{EPS}
\epsilon (|\nabla \phi({\bf r})|)=\epsilon_{m}+\frac{\epsilon_{s}-\epsilon_{m}}{(1+\alpha \frac{|\nabla \phi({\bf r})|^{2}}{2k_{B}T})^{p}}
\end{equation}
where the dielectric constant $\epsilon$ is piecewisely defined as $\epsilon=\epsilon_m$ in the macromolecule and $\epsilon=\epsilon_s$ in the solvent, $\alpha$ is a scaling parameter, and $p=1$ or 2 for different energy functions \cite{Hu12}. Another form is given as 
\begin{equation}\label{EPS2}
\epsilon (|\nabla \phi({\bf r})|)=\epsilon_{m}+(\epsilon_{s}-\epsilon_{m})
\exp( - \frac{|\nabla \phi({\bf r})|^{2}}{2k_{B}T} ).
\end{equation}
Mathematically, it can be verified that for both forms of the nonlinear dielectric function, we have
\begin{equation}
\lim_{| \nabla \phi | \to 0} \epsilon (|\nabla \phi|) = \epsilon_s, \quad \textrm{and}\quad
\lim_{| \nabla \phi | \to \infty} \epsilon (|\nabla \phi|) = \epsilon_m.
\end{equation}
Thus, for a finite solution, we have $\epsilon_m < \epsilon (|\nabla \phi|) < \epsilon_s$. In other words, the dielectric profile $\epsilon (|\nabla \phi|)$ undergoes a continuous transition from $\epsilon_m$ to $\epsilon_s$ at the solvent-solute interface when $|\nabla \phi|$ is continuous. 

To conduct numerical simulations, a finite domain $\Omega$ is usually chosen as the computational domain.  We may assign values along the boundary $\partial \Omega$ according to the approximate analytical condition
\begin{equation}\label{BC}
\phi ( {\bf r} ) = \frac{e^2_c}{k_B T} \displaystyle\sum\limits_{i=1}^{N_m} \frac{q_i}{\epsilon_s | {\bf r} - {\bf r}_i |} .
\end{equation}
When $\partial \Omega$ is of sufficient distance from the macromolecule, Eq. (\ref{BC}) can be utilized to approximate the results for potentials found from Eq. (\ref{NPE}). We note here that Eq. (\ref{BC}), for a collection of $N_m$ partial charges $q_i$ at positions ${\bf r}_i$, is simply a linear superposition of Coulomb's Law. For simplicity, the boundary $\partial \Omega$ is assumed to be of a cubic shape.

\subsection{Previous numerical approach}
A boundary value problem (BVP) approach has been constructed in \cite{Hu12} to numerically solve the BVP (\ref{NPE}) and (\ref{BC}), and is briefly reviewed here. Consider a uniform mesh partition of the computational domain $\Omega$. Without the loss of generality, we assume an equal grid spacing $h$ in all $x$, $y$ and $z$ directions. To facilitate the following discussions, we apply the notation $\phi_{i,j,k}=\phi(x_i,y_j,z_k)$ to denote the electrostatic potential at node $(x_i,y_j,z_k)$. 
In the BVP approach, a standard second order central finite differenc scheme is applied for the spatial discretization,  
\begin {align}\label{disnpe}
\epsilon(x_{i+\frac{1}{2}}, y_j, z_k)[ \phi_{i+1,j,k}-\phi_{i,j,k}]+\epsilon(x_{i-\frac{1}{2}}, y_j, z_k)[ \phi_{i-1,j,k}-\phi_{i,j,k}]\nonumber\\
+\epsilon(x_i, y_{j+\frac{1}{2}}, z_k)[ \phi_{i,j+1,k}-\phi_{i,j,k}]+\epsilon(x_i, y_{j-\frac{1}{2}}, z_k)[ \phi_{i,j-1,k}-\phi_{i,j,k}]\nonumber\\
+\epsilon(x_i, y_j, z_{k+\frac{1}{2}})[ \phi_{i,j,k+1}-\phi_{i,j,k}]+\epsilon(x_i, y_j, z_{k-\frac{1}{2}})[ \phi_{i,j,k-1}-\phi_{i,j,k}]&=-Q(x_i, y_j, z_k)h^{2}
\end {align}
where $Q(x_i, y_j, z_k)$ is the distribution of all source charges in the source term $\rho_m$, distributed by a trilinear interpolation. The permittivity $\epsilon$ on half grid nodes is approximated by an average over two nearby collocation nodes, e.g.,
\begin{equation}\label{EPSDS-0}
\epsilon(x_{i+\frac{1}{2}}, y_j, z_k)=
{\epsilon (|\nabla \phi|)}_{i+\frac{1}{2},j,k}=
\frac{1}{2} {\epsilon (|\nabla \phi|)}_{i,j,k} + 
\frac{1}{2} {\epsilon (|\nabla \phi|)}_{i+1,j,k},
\end{equation}
where ${\epsilon (|\nabla \phi|)}_{i,j,k}$ and ${\epsilon (|\nabla \phi|)}_{i+1,j,k}$ are approximated by the standard central differences.

In the BVP approach, one calculates the electrostatic potential $\phi$ and dielectric function $\epsilon$ alternatively. Normally, an initial dielectric profile which is a piecewise function based on $\epsilon_m$ and $\epsilon_s$ is first constructed. Then, one solves (\ref{disnpe}) by using a biconjugate gradient iterative method. With the calculated $\phi(x_i,y_j,z_k)$, one can update the dielectric function $\epsilon(x_i,y_j,z_k)$ and calculate half node values according to (\ref{EPSDS-0}). 
The computation will be regarded to be convergent, when the relative difference in the electrostatic solvation free energy of two successive steps is less than a tolerance. 
However, for extra large protein systems, the algebraic iterative solution process could be time consuming. 

\section {Alternating direction implicit (ADI) algorithm}
Motivated by our recent developments on the alternating direction implicit (ADI) algorithms for solving nonlinear Poisson-Boltzmann equations \cite{Zhao14,Geng13,Wilson16} and geometric flow equations in biomolecular surface generation \cite{Tian14}, we propose a pseudo-transient continuation approach in this paper to solve the nonlinear Poisson equation ({\ref{NPE}}). Mathematically, by introducing a pseduo-time, we construct a nonlinear diffussion equation
 \begin{equation}\label{STNPLE}
\frac{\partial \phi}{\partial t} ({\bf r},t)
=\nabla \cdot (\epsilon (|\nabla \phi({\bf r},t)|) \nabla \phi({\bf r},t))+\rho_m({\bf r}). 
\end{equation}
The solution to the original BVP (\ref{NPE}) and (\ref{BC}) can be obtained by a long time integration of (\ref{STNPLE}) until the steady state with boundary condition (\ref{BC}) and an appropriate initial condition. The nonlinearity can be effectively treated in such a time stepping process, provided that fully implicit or semi-implicit schemes are properly desigend \cite{Zhao14,Geng13,Tian14,Wilson16}.

In the present study, a Douglas-Rachford ADI scheme is proposed to solve the nonlinear heat equation (\ref{STNPLE}). Consider a uniform partition in time with an increment $\Delta t$. We will formulate the numerical scheme within the interval $[t_n, t_{n+1}]$ with $t_n=n\Delta t$. We denote $\phi^n_{i,j,k}=\phi(x_i,y_j,z_k,t_n)$. The nonhomogeneous diffusion process (\ref{STNPLE}) can be written explicitly in its Cartesian form
\begin{equation}\label{ADI3NPLE1}
\frac{\partial \phi}{\partial t}
=\frac{\partial}{\partial x}(\epsilon\frac{\partial \phi}{\partial x})
+\frac{\partial}{\partial y}(\epsilon\frac{\partial \phi}{\partial y})
+\frac{\partial}{\partial z}(\epsilon\frac{\partial \phi}{\partial z})
+\rho_m,
\end{equation}
The discretization of ({\ref{ADI3NPLE1}) using backward-Euler integration in time
and central differencing in space results in 
\begin{equation}\label{ADIEuler}
\phi_{i,j,k}^{n+1}=\phi_{i,j,k}^{n}+
\Delta t(\delta_{x}^{2}+\delta_{y}^{2}+\delta_{z}^{2})\phi_{i,j,k}^{n+1}+
\Delta tQ(x_{i},y_{j},z_{k}),
\end{equation}
where $\delta_{x}^{2}$, $\delta_{y}^{2}$ and $\delta_{z}^{2}$ are the central difference operators in the $x$, $y$ and $z$ directions, respectively. 
In particular, these operators are given as
\begin{equation}\label{phixx}
\left[\frac{\partial}{\partial x}(\epsilon\frac{\partial \phi}{\partial x})
\right]^{n+1}_{i,j,k} \approx \delta_{x}^2 \phi^{n+1}_{i,j,k}
:=\frac{1}{h^2}\Big( \epsilon^n_{i+\frac{1}{2},j,k}
(\phi^{n+1}_{i+1,j,k} - \phi^{n+1}_{i,j,k})
+\epsilon^n_{i-\frac{1}{2},j,k}
(\phi^{n+1}_{i-1,j,k} - \phi^{n+1}_{i,j,k}) \Big),
\end{equation}
\begin{equation}\label{phiyy}
\left[\frac{\partial}{\partial y}(\beta\frac{\partial \phi}{\partial y})
\right]^{n+1}_{i,j,k} \approx \delta_{y}^2 \phi^{n+1}_{i,j,k}
:=\frac{1}{h^2}\Big( \epsilon^n_{i,j+\frac{1}{2},k}
(\phi^{n+1}_{i,j+1,k} - \phi^{n+1}_{i,j,k})
+\epsilon^n_{i,j-\frac{1}{2},k}
(\phi^{n+1}_{i,j-1,k} - \phi^{n+1}_{i,j,k}) \Big),
\end{equation}
\begin{equation}\label{phizz}
\left[\frac{\partial}{\partial z}(\beta\frac{\partial \phi}{\partial z})
\right]^{n+1}_{i,j,k} \approx \delta_{z}^2 \phi^{n+1}_{i,j,k}
:=\frac{1}{h^2}\Big( \epsilon^n_{i,j,k+\frac{1}{2}}
(\phi^{n+1}_{i,j,k+1} - \phi^{n+1}_{i,j,k})
+\epsilon^n_{i,j,k-\frac{1}{2}}
(\phi^{n+1}_{i,j,k-1} - \phi^{n+1}_{i,j,k}) \Big).
\end{equation}
The three-dimensional (3D) system (\ref{ADIEuler}) can then be reduced to independent one-dimensional (1D) systems by the Douglas-Rachford ADI scheme
\begin{align}\label{ADIEulersplit}
\left(1-\Delta t \delta_{x}^{2}\right)\phi_{i,j,k}^{\star}&=
[1+\Delta t(\delta_{y}^{2}+\delta_{z}^{2})]\phi_{i,j,k}^{n}+\Delta t Q\left(x_{i}, y_{j}, z_{k}\right),\\
\left(1-\Delta t \delta_{y}^{2}\right)\phi_{i,j,k}^{\star\star}&=\phi_{i,j,k}^{\star}-\Delta t \delta_{y}^{2}\phi_{i,j,k}^{n},\\
\left(1-\Delta t \delta_{z}^{2}\right)\phi_{i,j,k}^{n+1}&=\phi_{i,j,k}^{\star\star}-\Delta t \delta_{z}^{2}\phi_{i,j,k}^{n}.
\end{align}

The key to the present ADI discretization is the treatment of $\epsilon$ function, which is the nonlinear term in the Poisson equation (\ref{NPE}). In time discretization, we note that $\epsilon$ functions in (\ref{phixx}), (\ref{phiyy}), and (\ref{phizz}) are evaluated at $t_n$, instead of $t_{n+1}$. This is because nonlinear algebraic systems will be resulted if evaluating at $t_{n+1}$, which introduces considerable difficulties. The approximation of $\epsilon^{n+1}$ by $\epsilon^n$ is of the order one, which is consistent with the temporal order of the implicit Euler scheme.

In spatial discretization, we will investigate two second order accurate approaches to evaluate $\epsilon$ at half nodes. Take (\ref{EPS}) as an example, and (\ref{EPS2}) will be similarly treated. In the first approach, we evaluate (\ref{EPS}) at half nodes, e.g.
\begin{equation}\label{epsilon1}
\epsilon^n_{i+\frac{1}{2},j,k}
 = \left[\epsilon_{m}+\frac{\epsilon_{s}-\epsilon_{m}}{(1+\alpha \frac{|\nabla \phi|^{2}}{2k_{B}T})^{p}}\right]^n_{i+\frac{1}{2},j,k},
\end{equation}
where$|\nabla \phi|^{2}=\phi_{x}^2+\phi_{y}^2+\phi_{z}^2$ and the subscripts of $\phi$ denote spatial derivatives, i.e., $\phi_x =\frac{\partial \phi}{\partial x}$. Following \cite{Tian14}, the discretization of  these derivatives are given by the following central difference formula
\begin{align}
\left(\phi_{x}^{2}\right)^n_{i+\frac{1}{2},j,k} & \approx
\left( \frac{\phi^{n}_{i+1,j,k} - \phi^{n}_{i,j,k}}{h} \right)^2, \nonumber \\
\left(\phi_{y}^{2}\right)^n_{i+\frac{1}{2},j,k} & \approx
\left( \frac{\phi^{n}_{i,j+1,k} - \phi^{n}_{i,j-1,k}}{4h}
+ \frac{\phi^{n}_{i+1,j+1,k} - \phi^{n}_{i+1,j-1,k}}{4h} \right)^2, \\
\left(\phi_{z}^{2}\right)^n_{i+\frac{1}{2},j,k} & \approx
\left( \frac{\phi^{n}_{i,j,k+1} - \phi^{n}_{i,j,k-1}}{4h}
+ \frac{\phi^{n}_{i+1,j,k+1} - \phi^{n}_{i+1,j,k-1}}{4h} \right)^2 \nonumber.
\end{align}
This approximation approach for $\epsilon$ will be noted as $\epsilon_{I}$ method in the following. 

In the second approach, $\epsilon^n_{i+\frac {1}{2},j,k}$ in (\ref {epsilon1}) will be first approximated via an average
\begin{align}
\epsilon^n_{i+\frac{1}{2},j,k}
 & =\frac{1}{2}\left( \epsilon_{m}+\frac{\epsilon_{s}-\epsilon_{m}}{(1+\alpha \frac{|\nabla \phi_{i,j,k}^{n}|^{2}}{2k_{B}T})^{p}}+\epsilon_{m}+\frac{\epsilon_{s}-\epsilon_{m}}{(1+\alpha \frac{|\nabla \phi_{i+1,j,k}^{n}|^{2}}{2k_{B}T})^{p}})\right),
\end{align}
where
\begin{align}
|\nabla \phi_{i,j,k}^{n}|^{2} & =((\phi_{x})_{i,j,k}^n)^{2}+((\phi_{y})_{i,j,k}^n)^{2}
+((\phi_{z})_{i,j,k}^n)^{2} \nonumber \\
& =\left( \frac{\phi^{n}_{i+1,j,k} - \phi^{n}_{i-1,j,k}}{2h} \right)^2+
\left( \frac{\phi^{n}_{i,j+1,k} - \phi^{n}_{i,j-1,k}}{2h} \right)^2+
\left( \frac{\phi^{n}_{i,j,k+1} - \phi^{n}_{i,j,k-1}}{2h} \right)^2,
\end{align}
and
\begin{align}
|\nabla \phi_{i+1,j,k}^{n}|^{2} & =((\phi_{x})_{i+1,j,k}^n)^{2}+((\phi_{y})_{i+1,j,k}^n)^{2}
+((\phi_{z})_{i+1,j,k}^n)^{2} \nonumber \\
& =\left( \frac{\phi^{n}_{i+2,j,k} - \phi^{n}_{i,j,k}}{2h} \right)^2+
\left( \frac{\phi^{n}_{i+1,j+1,k} - \phi^{n}_{i+1,j-1,k}}{2h} \right)^2+
\left( \frac{\phi^{n}_{i+1,j,k+1} - \phi^{n}_{i+1,j,k-1}}{2h} \right)^2.
\end{align}
We use notation $\epsilon_{II}$ for this approximation scheme. Note that the $\epsilon_{II}$ method is the one utilized in the original BVP approach \cite{Hu12}.

By eliminating $\phi_{i,j,k}^{\star}$ and $\phi_{i,j,k}^{\star\star}$ in ({\ref{ADIEulersplit}}) and solving for $\phi_{i,j,k}^{n+1}$, we obtain 
\begin{align}
\phi_{i,j,k}^{n+1}&=\phi_{i,j,k}^{n}+\Delta t(\delta_{x}^{2}+\delta_{y}^{2}+\delta_{z}^{2})\phi_{i,j,k}^{n+1}
+\Delta Q(x_{i},y_{j},z_{k}) \nonumber \\
&-\Delta t^{2}(\delta_{x}^{2}\delta_{y}^{2}+
\delta_{y}^{2}\delta_{z}^{2}+\delta_{x}^{2}\delta_{z}^{2})(\phi_{i,j,k}^{n+1}-\phi_{i,j,k}^{n})
+\Delta t^{3}\delta_{x}^{2}\delta_{y}^{2}\delta_{z}^{2}(\phi_{i,j,k}^{n+1}-\phi_{i,j,k}^{n})
\end{align}
In other words, the Douglas-Rachford scheme ({\ref{ADIEulersplit}}) is a higher order perturbation of the backward Euler scheme (\ref{ADIEuler}). Since both the backward Euler scheme and the temporal approximation of $\epsilon$ are first order in time, the proposed ADI scheme is of first order accuracy in time. The entire ADI time integration is semi-implicit, due to the evaluation of $\epsilon$ at $t_n$. Such an integration is unconditionally stable when the solution $\phi$ is sufficiently smooth. However, for real biological applications, the absolute stability may not be maintained. Computationally, each of 1D linear systems in ({\ref{ADIEulersplit}}) has a tridiagonal structure and can be efficiently solved by the Thomas algorithm. Both $\epsilon_{I}$ and $\epsilon_{II}$ methods are second order accurate in space. And the central finite difference approximations  in (\ref{phixx}), (\ref{phiyy}), and (\ref{phizz}) are also second order. Thus, the proposed ADI-$\epsilon_{I}$ and ADI-$\epsilon_{II}$ schemes have a second order  convergence in space and first order convergence in time.

\section{Numerical validation}
In this section, we validate the proposed ADI schemes by solving a nonlinear Poisson equation with a smooth analytical solution. We will explore the stability as well as spatial and temporal convergences of the ADI schemes. All computations are conducted on an SGI Xeon E5-4640 CPU core operating at 2.4 GHz and 8 GB of memory.

Consider a cubic domain $\Omega = [-\pi,\pi] \times[-\pi,\pi] \times[-\pi,\pi]$.
We solve a nonlinear Poisson equation of the form 
\begin{equation}\label{ADI3NPLE}
\frac{\partial \phi}{\partial t}
=\frac{\partial}{\partial x}(\epsilon\frac{\partial \phi}{\partial x})
+\frac{\partial}{\partial y}(\epsilon\frac{\partial \phi}{\partial y})
+\frac{\partial}{\partial z}(\epsilon\frac{\partial \phi}{\partial z})
+F,
\end{equation}
where a simplified form is adopted for the dielectric function
\begin{equation}
\epsilon (|\nabla \phi({\bf r})|)=\epsilon_{m}+\frac{\epsilon_{s}-\epsilon_{m}}
{1+\alpha |\nabla \phi({\bf r})|^{2}}.
\end{equation}
In the present study, the analytical solution is assumed to be
\begin{equation}\label{exacttest}
    \phi(x,y,z,t)=\sin(x)\sin(y)\sin(z)(1+e^{-\gamma t}).
\end{equation}
Thus, the source term $F$ is defined as
\begin{equation}\label{F}
    F(x,y,z,t)=\phi_{t}-\epsilon_{x}\phi_{x}-\epsilon\phi_{xx}
    -\epsilon_{y}\phi_{y}-\epsilon\phi_{yy}-
    \epsilon_{z}\phi_{z}-\epsilon\phi_{zz},
\end{equation}
where the derivatives $\phi_{t}$, $\phi_{x}$, $\phi_{xx}$, $\phi_{y}$,  $\phi_{z}$,  $\phi_{xy}$,  $\phi_{xz}$, 
 $\phi_{yy}$,  $\phi_{yz}$, $\phi_{zz}$ can be simply obtained by differentiating (\ref{exacttest}). In addition, the
$\epsilon_{x}$, $\epsilon_{y}$ and  $\epsilon_{z}$ are given by
\begin{equation}\label{epsilon_x}
\epsilon_{x}=(\epsilon_{m}-\epsilon_{s})\frac{2(\phi_{x}\phi_{xx}+\phi_{y}\phi_{yx}+\phi_{z}\phi_{zx})}{1+\phi_{x}^2+\phi_{y}^2+\phi_{z}^2},
\end{equation}
\begin{equation}\label{epsilon_y}
\epsilon_{y}=(\epsilon_{m}-\epsilon_{s})\frac{2(\phi_{x}\phi_{xy}+\phi_{y}\phi_{yy}+\phi_{z}\phi_{zy})}{1+\phi_{x}^2+\phi_{y}^2+\phi_{z}^2},
\end{equation}
\begin{equation}\label{epsilon_z}
\epsilon_{z}=(\epsilon_{m}-\epsilon_{s})\frac{2(\phi_{x}\phi_{xz}+\phi_{y}\phi_{yz}+\phi_{z}\phi_{zz})}{1+\phi_{x}^2+\phi_{y}^2+\phi_{z}^2}.
\end{equation}

\begin{table}[!tb]
\caption{Numerical convergence in space of the benchmark example.}
\label{tablenpespace}
\begin{center}
\begin{tabular}{llllllllllll}
\hline
 & \multicolumn{4}{c}{ADI-$\epsilon_I$} & \multicolumn{4}{c}{ADI-$\epsilon_{II}$}\\
\cline{2-5}\cline{6-9}
$h$ &$L_{\infty}$&Order&$L_{2}$&Order& $L_{\infty}$ & Order   & $L_{2}$ & Order \\
\hline
$\frac{\pi}{4}$  & 8.15E-02 &	     & 2.00E-02	&	     & 7.45E-02	&	     & 1.67E-02 &		\\
$\frac{\pi}{8}$  & 1.88E-02 & 2.12 & 5.50E-03	& 1.86 & 1.70E-02	& 2.13 & 4.23E-03	& 1.98\\
$\frac{\pi}{16}$ & 4.65E-03 & 2.02 & 1.45E-03	& 1.92 & 4.17E-03	& 2.03 & 1.10E-03	& 1.94\\
$\frac{\pi}{32}$ & 1.18E-03 & 1.98 & 3.80E-04 & 1.93 & 1.05E-03	& 1.99 & 2.87E-04	& 1.94\\
\hline
\end{tabular}
\end{center}
\end{table}

In our computations, we chose $\epsilon_s=80$, $\epsilon_m=1$, $\gamma=0.1$, and $\alpha=0.1$. The initial value of $\phi$ is chosen as the analytical solution at time $t=0$. The stopping time is set to be $T=10$. With the given analytical solution, the numerical  errors in $L_{\infty}$ and $L_{2}$ norms will be reported.

In the present study, the proposed ADI scheme is found to be unconditionally stable, since the solution is a smooth function. We thus will report only spatial and temporal accuracies. 
We first test the spatial order by using a sufficiently small $\Delta t=0.001$. 
The numerical errors with different $h$ are reported in Table \ref{tablenpespace}. It can be seen that both $\epsilon_I$ and $\epsilon_{II}$ methods achieve second order in space, and the difference between two methods is minor in this example.

\begin{table}[!tb]
\caption{Numerical convergence in time of the benchmark example.}
\label{tablenpetime}
\begin{center}
\begin{tabular}{llllllllllll}
\hline
 & \multicolumn{4}{c}{ADI-$\epsilon_{I}$} & \multicolumn{4}{c}{ADI-$\epsilon_{II}$}\\
\cline{2-5}\cline{6-9}
$\Delta t$ &$L_{\infty}$&Order&$L_{2}$&Order& $L_{\infty}$ & Order   & $L_{2}$ & Order \\
\hline
$0.8$ &6.43E-01 &	     & 2.22E-01	&      & 6.42E-01	&	     & 2.24E-01 &		\\
$0.4$ &6.05E-01 & 0.09 & 2.11E-01	& 0.07 & 6.04E-01	& 0.09 & 2.11E-01	& 0.09 \\
$0.2$ &4.65E-01 & 0.38 & 1.64E-01	& 0.36 & 4.65E-01	& 0.38 & 1.65E-01	& 0.35\\
$0.1$ &1.37E-01 & 1.76 & 5.05E-02 & 1.70 & 1.37E-01	& 1.76 & 5.03E-02	& 1.71\\
$0.05 $&2.33E-02& 2.56 & 7.22E-03 & 2.81 & 2.33E-02	& 2.56 & 7.16E-03	& 2.81\\
\hline
\end{tabular}
\end{center}
\end{table}

\begin{figure*}[!tb]
\begin{center}
\begin{tabular}{cc}
\psfig{figure=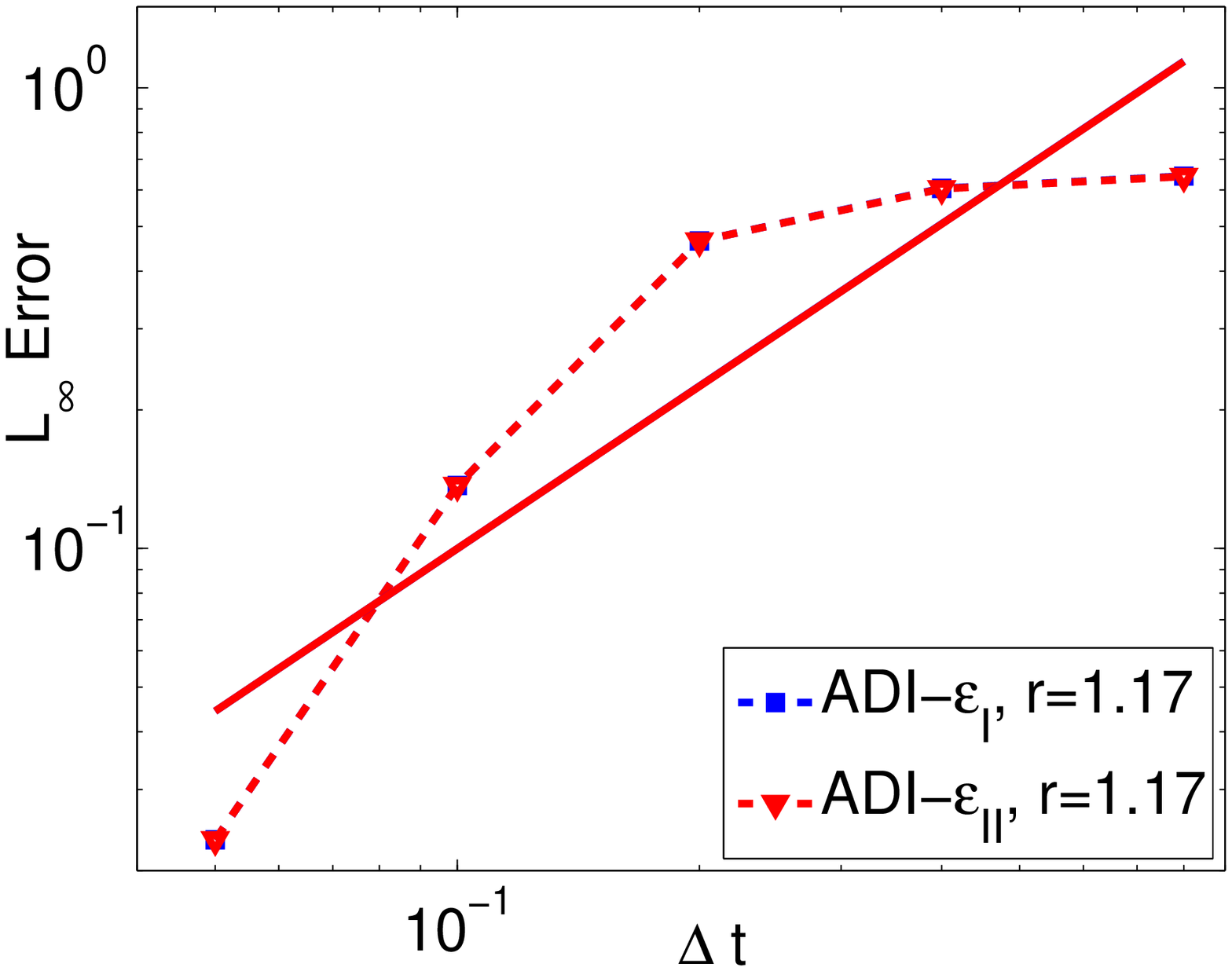,width=0.45\linewidth} &
\psfig{figure=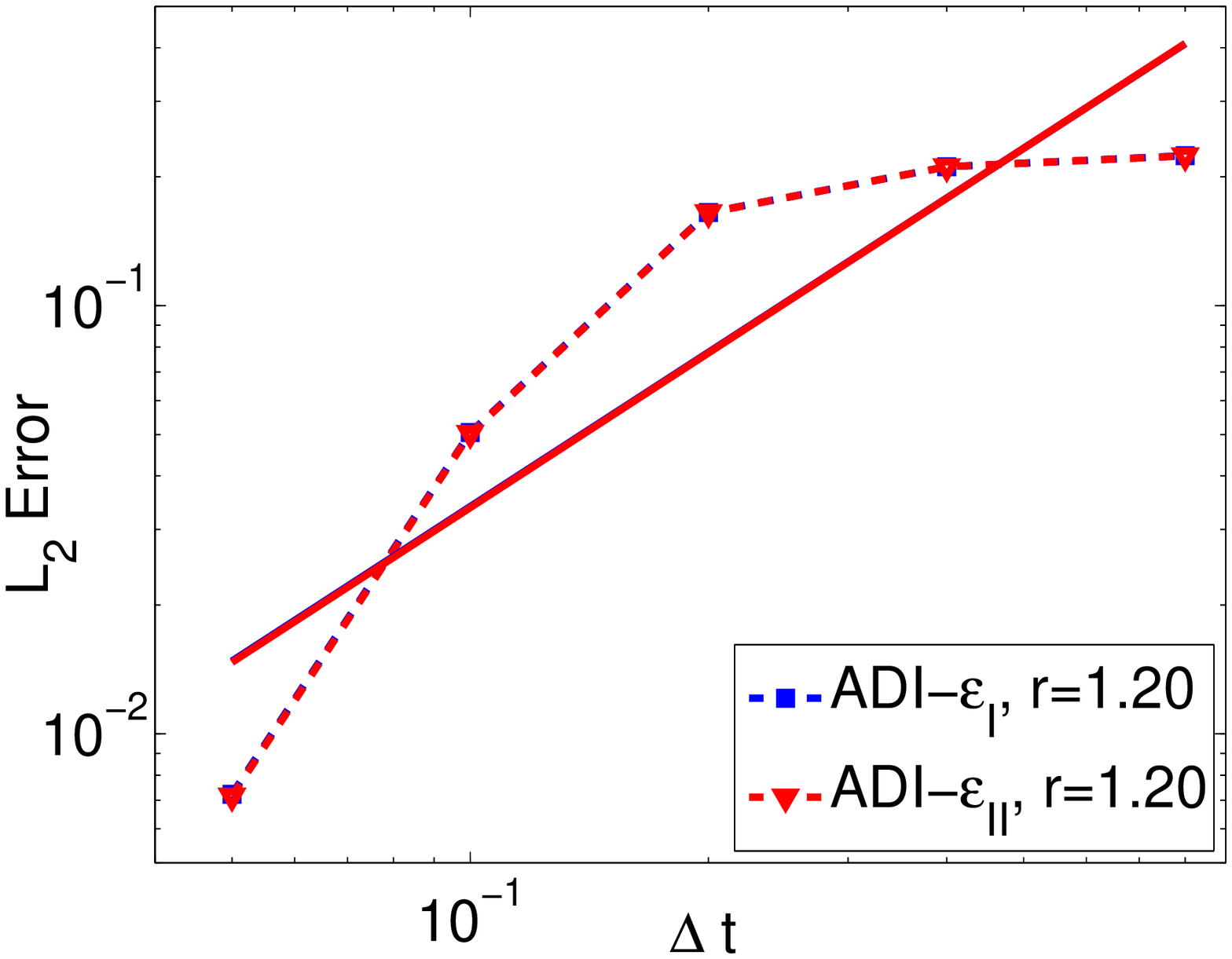,width=0.45\linewidth} \\
(a) & (b)
\end{tabular}
\end{center}
\caption{Temporal convergence tests of the benchmark example.
(a). $L_\infty$ error; (b). $L_2$ error.}
\label{fig.ex1}
\end{figure*}

In Table \ref{tablenpetime}, with $h=\frac{\pi}{48}$, numerical errors and temporal orders of the ADI schemes are reported. It can be observed that the convergence is quite slow when $\Delta t$ is large and becomes faster when $\Delta t$ decreases. Since such temporal orders are rather irregular, we further plot the $L_\infty$ and $L_2$ errors in Fig. \ref{fig.ex1}. In both charts, a least squares fitting is conducted in the log-log scale. The fitted convergence lines are shown as solid lines. The slopes $r$ of the solid straight lines are calculated and listed in the legends of both plots. The magnitude of $r$ reflects the average convergence rate. Fig. \ref{fig.ex1} indicates that for both $L_\infty$ and $L_2$ errors and both $\epsilon_I$ and $\epsilon_{II}$ methods, the overall order in time is about one. Since

\section{Application for solvation analysis}
In general, the electrostatic free energy of solvation is complemented by the nonpolar contribution. The nonpolar solvation free energy functional proposed by Wagoner and Baker \cite{Wagoner} is considered for the nonpolar contribution,
\begin{align}
G_{np}=\gamma(Area)+p(Vol)+\rho_0\int_{\Omega_s}\rho_{s}U^{att} d{\bf r}
\end{align} \label{Nonpolare1}
where $\gamma$ is the surface tension,
$p$ is the hydrodynamic pressure,
Vol is the volume occupied by the molecule,
$\rho_{s}$ is the solvent bulk density,
$\Omega_s$ denotes the solvent-accessible region,
and $U^{att}$ is the attractive portion of the van der Waals potential at {\bf {r}}.

With the given form of the nonlinear dielectric function, there is no sharp interface between the solvent and solute, which avoids the surface singularity and numerical instability. A linear mapping is introduced for the surface $S({\bf r})$\cite{Hu12}.
\begin{align}
S(\epsilon)=\frac{1}{\epsilon_{s}-\epsilon_{m}}(\epsilon_{s}-\epsilon(|\nabla \phi ({\bf r})|))
\end{align} \label{MappingS}
where $S(\epsilon)$ is a characteristic function such that the solute domain $\Omega_m$ is defined by $S(\epsilon)\ge0$. Similarly, 1-$S(\epsilon)$ is a solvent characteristic function such that the solvent domain $\Omega_s$ is defined as $1-S(\epsilon)\ge0$. Once the surface profile $S$ is determined, then the volume of solute can be expressed in terms of
surface $S({\bf r})$\cite{Hu12}.
\begin{align}
Vol=\int_{\Omega_m}d{\bf r}=\int_{\Omega}S(\epsilon)d{\bf r}
\end{align} \label{MappingV}
where $\Omega$ is the total domain of the solvation.
By using the coarea formula of geometric measure theory \cite{Wei10}, the surface area can also be expressed
as a volume integral
\begin{align}
Area=\int_{\Omega}|\nabla S(\epsilon)|d{\bf r}
\end{align} \label{Mappingarea}
It is noted that $Area$ only has  contributions from a transition region of the solvent-solute boundary.  The van der Waals dispersion term can be rewritten as
\begin{align}
\int_{\Omega_s}\rho_{s}U^{att}d{\bf r}=\int_{\Omega}\rho_{s}(1-S(\epsilon))U^{att}d{\bf r}
\end{align} \label{Mappingarea2} 
Therefore, a practical way to compute nonpolar free energy $G_{np}$ is to express in terms of the as surface 
profile $S(\epsilon)$\cite{Hu12},
\begin{align}
G_{np}=\int_{\Omega} \gamma|\nabla S(\epsilon)|+pS(\epsilon)+\rho_s(1-S(\epsilon))U^{att}(\bf r) d{\bf r}
\end{align} \label{Nonpolarec}

With the steady state values of
the hypersurface function $S$ and the electrostatic potential $\phi$,
the solvation free energy is calculated as the follows.
It is known that the total free energy functional of solvation does not directly provide
the total solvation free energy.
Actually, one needs to calculate the difference of the
macromolecular system in the vacuum and in the solvent.
The solvation free energy can be computed as
\begin{equation}\label{totalDG}
\Delta G = G_{\rm np} + ( G_{\rm p} - G_0).
\end{equation}
where $G_{\rm np}$ and $G_{\rm p}$ are, respectively, the nonpolar and polar
solvation free energies of the solute solvent system with different $\epsilon_s$
and $\epsilon_m$, while $G_0$ is the polar free energy calculated from the
homogeneous  (vacuum) environment with $\epsilon_s=\epsilon_m=1$ \cite{Zhao11}.
The term
$G_{\rm p} -G_0$ can be regarded as the electrostatic solvation free energy.
The nonpolar solvation free energy of the macromolecule, $G_{\rm np}$,
is computed exactly according to its definition \cite{Zhao11}.
In the present study, the polar part is evaluated as
\begin{equation}
G_{\rm p} = \frac{1}{2} \int_{\Omega} S({\bf r}) \rho_m \phi({\bf r}) d{\bf r}
= \frac{1}{2} \sum_{i=1}^{N_m} Q({\bf r}_i) \phi ({\bf r}_i),
\end{equation}
where $Q({\bf r}_i)$ is the $i$th partial charge at location ${\bf r}_i$ in the
biomolecule, and $N_m$ is the total number of partial charges.
Similarly, the electrostatic solvation free energy can be calculated as
\begin{equation}
\Delta G_{\rm p} = G_{\rm p} - G_0 = \frac{1}{2} \sum_{i=1}^{N_m} Q({\bf r}_i) (\phi ({\bf r}_i) - \phi_0({\bf r}_i) ),
\end{equation}
where $\phi$ and $\phi_0$ are electrostatic potentials in the presence of the
solvent and the vacuum, respectively.

We validate our ADI scheme through  one atom of unit van der Waals radius and unit charge in the solvent. Both the atomic center and 
the charge are located at the origin point. We choose $\epsilon_{m}$=1 and $\epsilon_{s}$=80 for dielectric constants in the 
solute and the solvent, respectively. A fine mesh size of $h=0.25$ \AA ~ is used in our computation. In our benchmark experiment and our one atom test, both $\epsilon_I$ and $\epsilon_{II}$ almost have identical results, therefore, we will only report the results from $\epsilon_I$  in the following sections.  Due to a nonlinear dielectric function in the NPE system and nonsmooth solution in molecular system, we will test the stability of ADI scheme by applying it into one unit atom system.  

 \begin{figure*}[!tb]
\begin{center}
\begin{tabular}{cc}
\psfig{figure=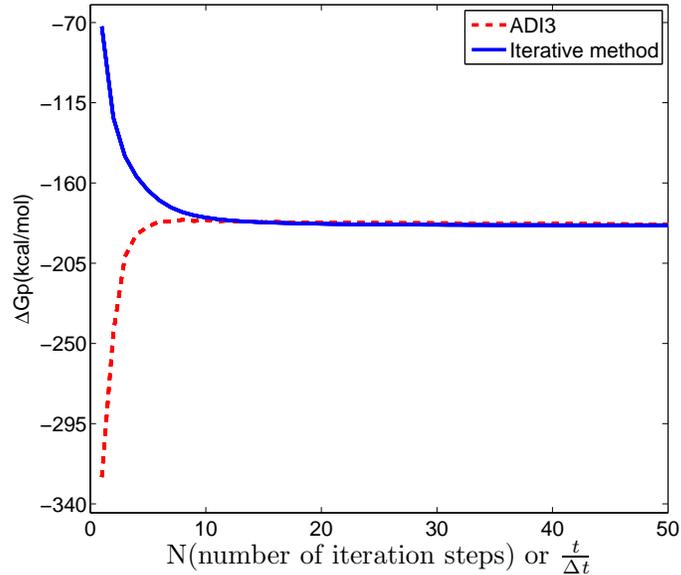,width=0.6\linewidth} 
\end{tabular}
\end{center}
\caption{The convergence of electrostatic solvation free energies $\Delta G_{p}$ over the iteration steps $N$ or iteration time $t$ of one-atom system.}
\label{fig.areaoneatom}
\end{figure*}
  \begin{figure*}[!tb]
\begin{center}
\begin{tabular}{cc}
\psfig{figure=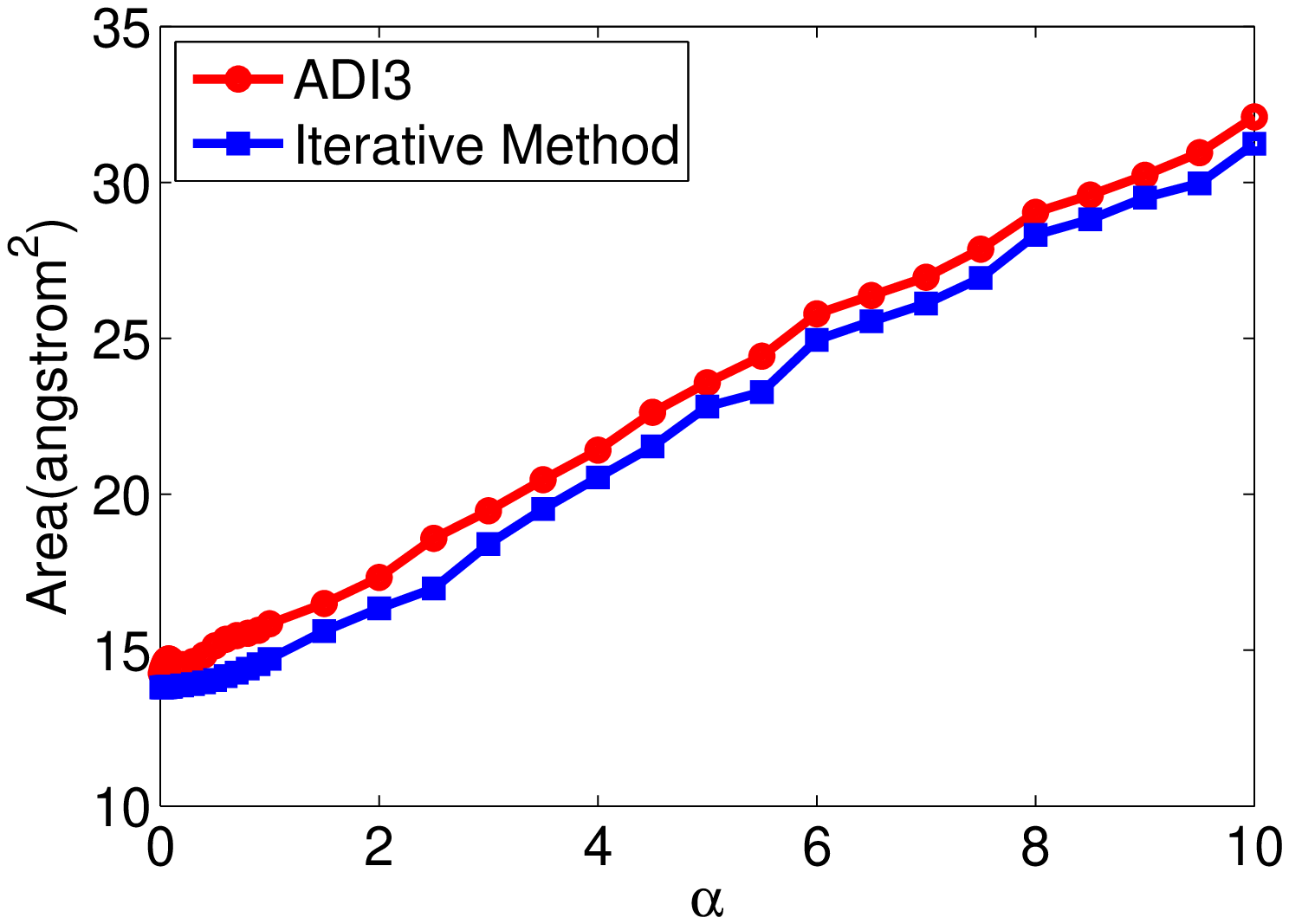,width=0.45\linewidth} &
\psfig{figure=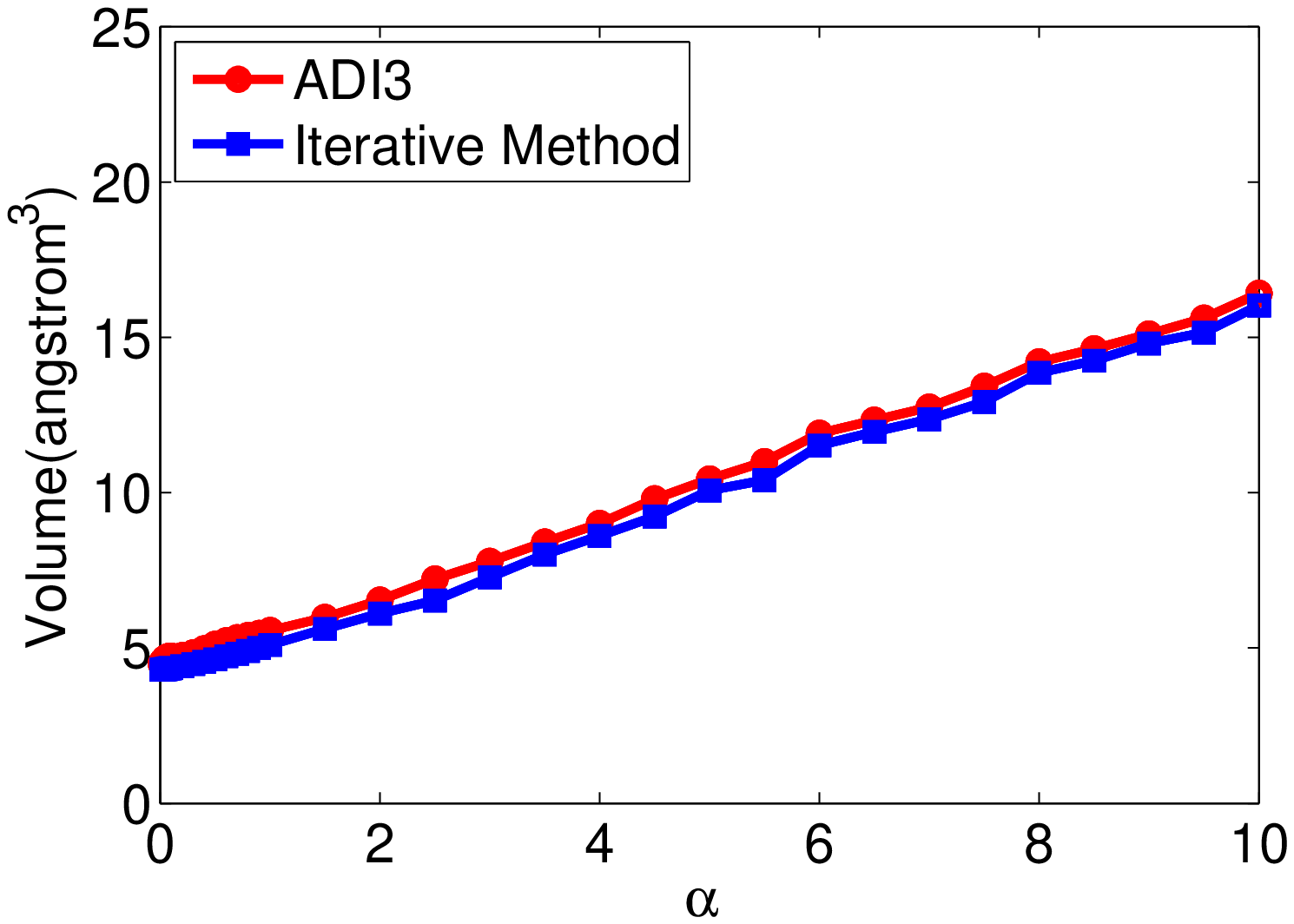,width=0.45\linewidth} \\
(a) & (b)
\end{tabular}
\end{center}
\caption{Area and volume of the nonliear Poisson model solved by the ADI and iterative methods for a wide range of $\alpha$ values of a one-atom system
with atomic radius 1$\AA$, $\epsilon_{m}$=1 and $\epsilon_{s}=80$:
(a) Comparison of area;
(b) Comparison of Volume.}
\label{fig.volumeoneatom}
\end{figure*}

From  Figures \ref{fig.areaoneatom} and \ref{fig.volumeoneatom}, we can see that the iterative method with $N$=50 steps and the ADI method with T=5 and $\Delta t$=0.1 converge to almost the same electrostatic solvation energy. Since $\alpha$ is a parameter to adjust the strength of hyperpolarizability, it is worthwhile to further explore the impact of  hyperpolarization by testing different $\alpha$ values. By setting up the tolerance $\tau$=0.01 (kcal/mol) for iterative method and $T=2$ and $\Delta t=0.05$ for ADI scheme, which has the same level of successive iteration difference. We next validate the area, volume and  electrostatic solvation energy of one atom solved by both schemes with a wide range of $\alpha$ values.

  \begin{figure*}[!tb]
\begin{center}
\begin{tabular}{cc}
\psfig{figure=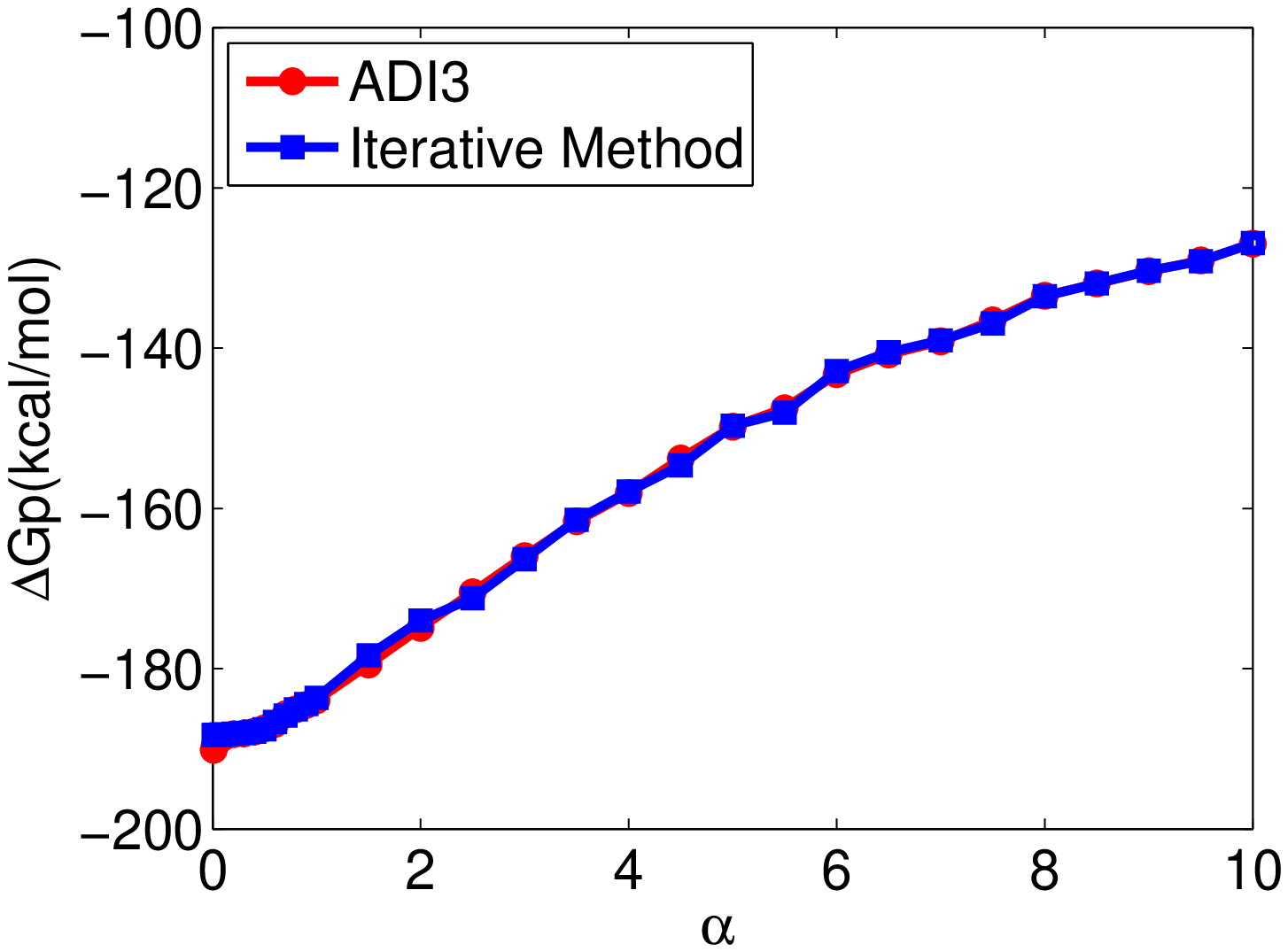,width=0.45\linewidth} &
\psfig{figure=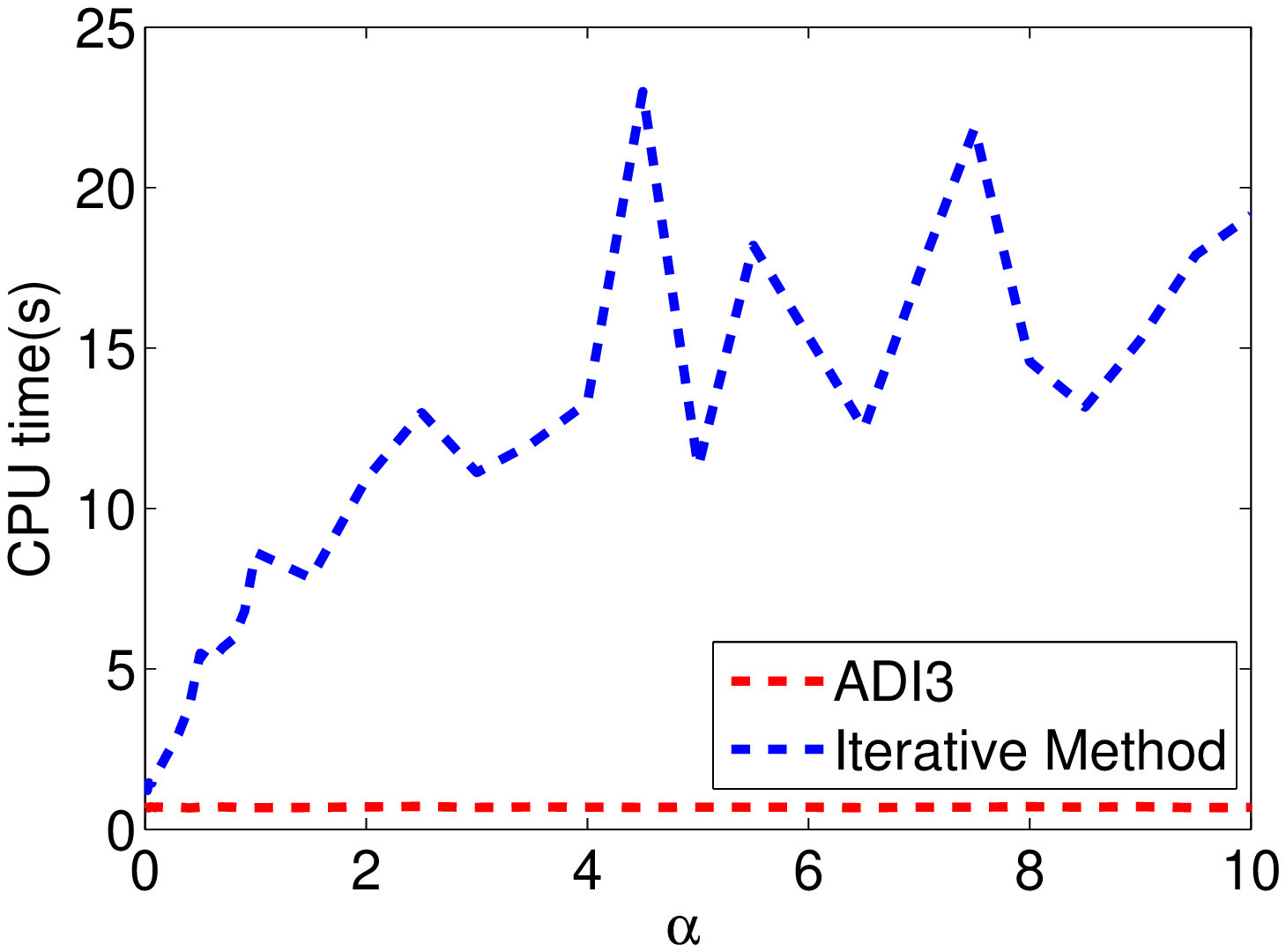,width=0.45\linewidth} \\
(a) & (b)
\end{tabular}
\end{center}
\caption{Electrostatic solvation energy of the nonlinear Poisson model and CPU time by the
ADI scheme and the iterative method for a wide range of $\alpha$ values for a one-atom system with atomic 
radius 1 \AA, $\epsilon_{m}$=1 and $\epsilon_{s}=80$:
(a) Electrostatic solvation energy ;
(b) CPU time.}
\label{fig.ADIiterativecompare}
\end{figure*}

\begin{table}[!tb]
\caption{Electrostatic solvation free energies (kcal/mol) for 17 compounds.}
\label{table.17setAI}
\begin{center}
\begin{tabular}{lrrrrr}
\hline
 & &  \multicolumn{2}{c}{BCG-NPE\cite{Hu12}}  &  \multicolumn{2}{c}{Present} \\
\cline{3-4} \cline{5-6}
Compound & Exptl &  $\Delta G$ & Error &  $\Delta G$ & Error   \\
\hline
glycerol triacetate & $-8.84$  & $-10.42$ & $-1.58$ & $-9.90$ & $-1.06$ \\
benzyl bromide & $-2.38$ & $-3.47$ & $-1.09$ & $-3.98$ & $-1.60$ \\
benzyl chloride & $-1.93$ & $-3.65$ & $-1.72$ & $-4.13$ & $-2.20$ \\
m-bis(trifluoromethyl)benzene &  $1.07$  & $-0.98$ & $-2.05$ & $-0.94$ & $-2.01$ \\
N,N-dimethyl-p-methoxybenzamide & $-11.01$ & $-7.26$ & $3.75$
& $-7.27$ & $3.74$ \\
N,N-4-trimethylbenzamide & $-9.76$ & $-5.81$ & $3.95$
& $-6.01$ & $3.75$\\
bis-2-chloroethyl ether & $-4.23$ &  $-2.68$ & $1.55$
& $-2.67$ & $1.57$ \\
1,1-diacetoxyethane & $-4.97$ &  $-6.63$ & $-1.66$
& $-6.55$ & $-1.58$\\
1,1-diethoxyethane & $-3.28$  & $-2.93$ & $0.35$
& $-2.84$ & $0.44$\\
1,4-dioxane & $-5.05$ &  $-4.58$ & $0.47$
& $-4.48$ & $0.57$ \\
diethyl propanedioate & $-6.00$ & $-6.14$ & $-0.14 $
& $-6.13$ & $-0.13$\\
dimethoxymethane & $-2.93$ & $-3.53$ & $-0.60$
& $-3.53$ & $-0.60$ \\
ethylene glycol diacetate & $-6.34$  & $-7.04$ & $-0.70 $
& $-6.82$ & $-0.48$  \\
1,2-diethoxyethane & $-3.54$  & $-2.64$ & $0.90$
& $-2.56$ & $0.98$ \\
diethyl sulfide & $-1.43$ &  $-1.06$ & $0.37$
& $-1.50$ & $0.07$\\
phenyl formate & $-4.08$  & $-6.51$ & $-2.43 $
& $-6.77$ & $-2.69$\\
imidazole & $-9.81$  & $-9.68$ &$0.13$
& $-9.52$ & $0.29$\\
\hline
RMSE (kcal/mol) &  &  \multicolumn{2}{c}{$1.7774$} &\multicolumn{2}{c}{$1.7985$} \\
 Average error (kcal/mol) &  &  \multicolumn{2}{c}{$1.3777$} &\multicolumn{2}{c}{$1.3979$}\\
\hline
\end{tabular}
\end{center}
\end{table}

To further validate the ADI scheme for the nonlinear Poisson model, we apply it to the solvation analysis of a set of 17 small compounds. Various approaches, including quantum mechanics and Poisson-Boltzmann theory, have been applied to this test set by Nicholls \emph{et al}. \cite{Nicholls}. This test set has also been employed in the earlier work to validate the ADI scheme on solving differential geometry-based solvation models\cite{Tian14, Zhao14, Zhao11}. This test set has been frequently used because its experimental data of solvation energies are available. Both polar and nonpolar solvation energies are required in order to compare with experimental data. A dense mesh with $h$=0.25 \AA~ is employed in order to achieve a better spatial resolution for these small chemical compounds. As we saw in the one-atom system, the solvation energy is converged at T=1.0. Therefore, we set  T=2.0 and $\Delta t=0.1$ for the ADI scheme.  $\alpha$=40 is chosen for both iterative method and the ADI scheme. The calculated solvation free energies are listed in Table \ref{table.17setAI}. The root mean-square error (RMS) of the computation results is 1.80 and the average error is 1.40, while the RMSE is 1.78 and the average error is 1.38 with the iterative scheme \cite{Hu12}. The existing explicit solvent approach, which is more expensive, reduces the RMSE to 1.71 kcal/mol \cite{Nicholls}. With both the iterative method and the ADI scheme, the nonlinear Poisson model under discussion provides a relatively good prediction of solvation energies for this set of molecules. However,  we can see from Fig. (\ref{fig.ADIiterativecompare2}), that the CPU time cost of the ADI scheme is much smaller than the iterative method. 
\begin{table}[!tb]
\caption{Electrostatic solvation free energies (kcal/mol) for 19 proteins.}
\label{table.19protein}
\begin{center}
\begin{tabular}{llllll}
\hline
PDB ID & No. of atoms & DGSM\cite{ChenZhan1} & BCG-NPE\cite{Hu12}
&ADI &  \\
\hline
1ajj & 519 & $-1178.5$ & $-1088.6$ & $-1211.2$&\\
2erl & 573 & $-935.8$  & $-957.0$  &  $-905.5$\\
1bbl & 576 & $-965.9$  & $-963.5$  & $-876.4$ \\
1vii & 596 & $-892$  & $-885.4$  &  $-805.6$\\
2pde & 667 & $-843$  & $-823.6$  & $-953.6$\\
1sh1 & 702 & $-771$  & $-779.2$  &  $-725.7$\\
1fca & 729 & $-1200.6$ & $-1252.0$ &  $-1197.0.0$\\
1uxc & 809 & $-1125.7$ & $-1119.0$ &  $-1021.2$\\
1fxd & 824 & $-3291.9$ & $-3343.7$ &  $-3262.6$\\
1bor & 832 & $-871.4$  & $-875.3$  & $-856.8$\\
1hpt & 858 & $-808.2$  & $-802.7$  &  $-718.8$\\
1mbg & 903 & $-1328.2$ & $-1338.4$ & $-1259.0$\\
1r69 & 997 & $-1072.7$ & $-1075.9$ & $-971.3$\\
1neq & 1187& $-1713.9$ & $-1710.8$ & $-1627.4$\\
451c & 1216& $-102406$ & $-1036.6$ &  $-1011.3$\\
1a2s & 1272& $-1900.3$ & $-1967.9$ &  $-1922.2$\\
1svr & 1435& $-1704.6$ & $-1959.5$ & $-1715.0$\\
1a63 & 2065& $-2380.5$ & $-2371.6$ &  $-2374.8$\\
1a7m & 2809& $-2179.8$ & $-2283.7$ &  $-2287.9$\\
\hline
\end{tabular}
\end{center}
\end{table}

  \begin{figure*}[!tb]
\begin{center}
\begin{tabular}{cc}
\psfig{figure=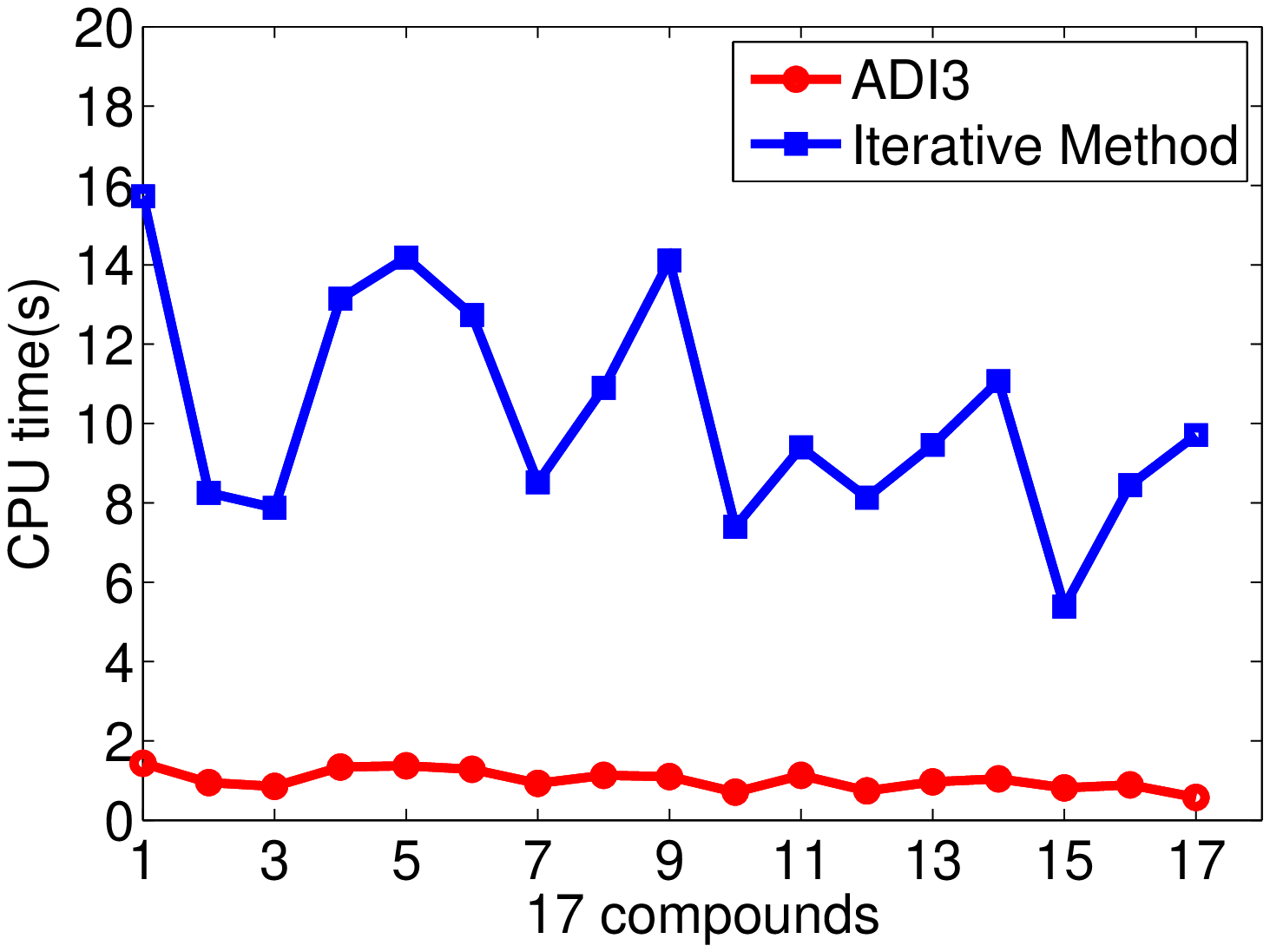,width=0.45\linewidth} & \quad
\psfig{figure=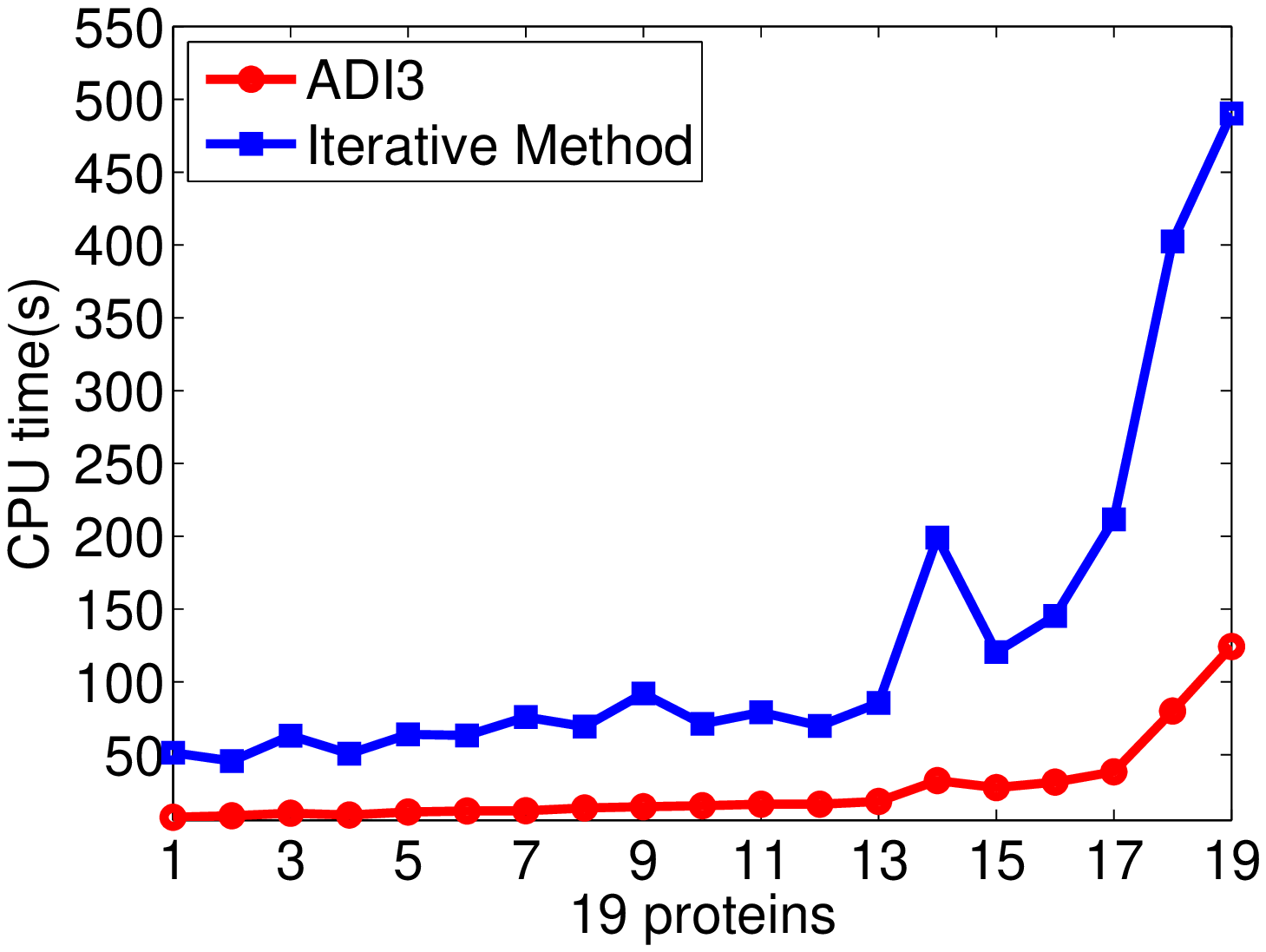,width=0.45\linewidth}\\
(a) &(b)
\end{tabular}
\end{center}
\caption{CPU time of the iterative method and the ADI scheme for 17 compounds and 19 proteins with $\epsilon_{m}$=1 and $\epsilon_{s}=80.$}
\label{fig.ADIiterativecompare2}
\end{figure*}

After demonstrating the success in validation and application of solvation analysis in one atom and small compounds,  
we finally calculate the electrostatic solvation energy for the set of 19 proteins which have been studied in the previous
 chapters. The number of atoms for this set of proteins ranges from 519 to 2809. Here, a large spacing $h=0.5\AA$ is
 employed due the large size of the molecules. The ADI scheme parameters are set to be $\Delta t$=0.15 and $T$=3 but not $T=10$ when we set up successive iterations differ by less than 1 (kcal/mol). 
 Under the same level of convergence error 1 kcal/mol, the results for 19 proteins are summarized in Table
 \ref{table.19protein}. Also, we compare the electrostatic solvation free energies with the results of DGSM\cite{ChenZhan1} and the
 iterative method BCG-NPE \cite{Hu12}. It can be seen from Table \ref{table.19protein} that there is general agreement
 between the present results and the existing literature. For a fixed $\alpha$, the CPU acceleration is moderate. We can see from Fig. \ref{fig.ADIiterativecompare2} that our ADI scheme is also 4-5 times faster
 than the BCG interative method. 

\section{Conclusion}

We have applied ADI scheme in solving a nonlinear Poisson equation for electrostatic analysis. It is first validated 
by comparing it with analytical solutions. Using two different ways of discretization the dielectric function $\epsilon$, the proposed ADI scheme can be shown to be a second order of accuracy in space and about $1.5^{th}$ order in time. 
The ADI scheme is extensively validated with the Kirkwood model (one-atom) and experimental measurements of 17 small chemical compounds. There is a solid agreement between the ADI scheme and an existing iterative method on the solvation analysis of small molecules. Application of the proposed scheme for a nonlinear Poisson model is also considered for electrostatic analysis of 19 proteins. Our ADI scheme is 4-5 times faster than existing BCG-iterative method on the application of a nonlinear Poisson equation for heterogeneous media.
 
The development of an efficient nonlinear solvation simulator based on the present
ADI molecular surface construction and a fast nonlinear Poisson-Boltzmann solver
is currently under our consideration.
Stability analysis of the entire solvation system and adaptive steady state
solution will be explored in the future.


\begin{thebibliography}{99}
\bibitem{Baker05} N. A. Baker,
Improving implicit solvent simulations: a Poisson-centric view,
\emph{Current Opinion in Structural Biology}, {\bf 15}, 137-143, (2005).

\bibitem{ChenZhan1} Z. Chen, N. Baker, and G.W. Wei,
Differential geometry based solvation model I: Eulerian formation,
\emph{J.  Comput. Phys.}, {\bf 229}, 8231-8258, (2010).

\bibitem{ChenZhan2} Z. Chen, N. Baker, and G.W. Wei,
Differential geometry based solvation model II: Lagrangian formation,
\emph{J.  Math. Bio.}, {\bf 63}, 1139-1200, (2011).

\bibitem{Crowley} P.B. Crowley and A. Golovin,
Cation-pi interactions in protein-protein interfaces,
\emph{Proteins - Struct. Func. Bioinf.},  {\bf 59}, 231-239 (2005).

\bibitem{Deng15} W. Deng, X. Zhufu, J. Xu and S. Zhao, A new discontinuous Galerkin
method for the nonlinear Poisson-Boltzmann equation A. \emph{Math. Lett.}, 257,
1000-1021, (2015).

\bibitem{Dong08} F. Dong, B. Olsen, and N.A. Baker,
Computational methods for biomolecular electrostatics,
\emph{Methods in Cell Biology}, {\bf 84}, 843-870, (2008).

\bibitem{Feig04} M. Feig, C. L. Brooks III,
Recent advances in the development and application of implicit solvent models
in biomolecule simulations,
\emph{Current Opinion in Structural Biology}, {\bf 14}, 217-224, (2004).

\bibitem{Geng13} W. Geng and S. Zhao,
Fully implicit ADI schemes for solving the nonlinear Poisson-Boltzmann equation,
\emph{Molecular Based Math. Bio.}, {\bf 1}, 109-123, (2013).

\bibitem{Hu12} L.Hu and G.Wei, Nonlinear Poisson equation for heterogeneous media, 
\emph{Biophysical Journal},{\bf 103}, 758-766,(2012).



\bibitem{Nicholls} A. Nicholls, D.L. Mobley, P.J. Guthrie, J.D. Chodera, and V.S. Pande,
Predicting small-molecule solvation free energies: An informal blind test for
computational chemistry,
\emph{J. Medicinal Chem.},   {\bf 51}, 769-779, (2008).

\bibitem{Roux} B. Roux, T.Simpson, 
Implicit solvent models,
\emph{Biophys. Chem} {\bf 78}, 1-20,(1999).



\bibitem{Sharp} Sharp KA, Honig B, 
Calculating total electrostatic energies with the nonlinear Poisson-Boltzmann equation. 
\emph{Journal of Physical Chemistry},{\bf 94}7684-7692, (1990).

\bibitem {Tian14} S. Zhao,
A fast ADI algorithm for geometric flow equations in biomolecular surface generation,
\emph{Int. J. Numer. Meth. Biomed. Engng}, {\bf 30}, 490-516, (2014).

\bibitem{Wei10} G.W. Wei,
Differential geometry based multiscale models,
\emph{Bullet. Math. Bio.}, {\bf 72}, 1562-1622, (2010).

\bibitem{Wilson16} L. Wilson and S. Zhao, Unconditionally stable time splitting methods for
the electrostatic analysis of solvated biomolecules, Int. J. Numer. Anal.
Model., 13, 852-878, (2016).


\bibitem{Wagoner} J. Wagoner and N.A. Baker,
Assessing implicit models for nonpolar mean solvation forces: the importance
of dispersion and volume terms,
\emph{Proceedings of the National Academy of Sciences of the USA},
{\bf 103}, 8331-8336, (2006).

\bibitem {Zhao11} S. Zhao,
Pseudo-time-coupled nonlinear models for biomolecular surface representation
and solvation analysis,
\emph{Int. J. Numer. Meth. Biomed. Engng}, {\bf 27}, 1964-1981, (2011).

\bibitem {Zhao14} S. Zhao,
Operator splitting ADI schemes for
pseudo-time coupled nonlinear solvation simulations,
\emph{J. Comput. Phys.}, {\bf 257}, 1000-1021, (2014).



\end{thebibliography}
\end{document}